\documentclass[12pt]{article}

\usepackage{graphicx, amssymb, latexsym, amsfonts, amsmath, lscape, amscd,
	amsthm, color, epsfig, mathrsfs, tikz, enumerate}
\usepackage{subcaption, geometry}
\usepackage{csquotes}

\usepackage{amssymb}
\usepackage{color}
\usepackage{mathtools}
\usepackage{pgf,tikz}
\usetikzlibrary{decorations.pathreplacing}
\usetikzlibrary{decorations.markings}
\usetikzlibrary{positioning}
\usetikzlibrary{arrows}

\newtheorem{theorem}{Theorem}[section]
\newtheorem{conjecture}[theorem]{Conjecture}

\newtheorem{lemma}[theorem]{Lemma}
\newtheorem{proposition}[theorem]{Proposition}
\newtheorem{question}[theorem]{Question}

\theoremstyle{definition}  

\newcommand\DELETE[1]{}

\begin{document}
	
	%\begin{frontmatter}
	
	\title{{\bf On coloring parameters of triangle-free planar $(n,m)$-graphs}}
	\author{
		{\sc Soumen Nandi}$\,^{a}$, {\sc Sagnik Sen}$\,^{b}$, {\sc S Taruni} $\,^{b}$ \\
		\mbox{}\\
		{\small $(a)$ Netaji Subhas Open University, India.}\\
		{\small $(b)$ Indian Institute of Technology Dharwad, India.}
	}

	\date{}
	
	\maketitle

	\begin{abstract}
		An $(n,m)$-graph is a graph with $n$ types of arcs and $m$ types of edges. A homomorphism of an $(n,m)$-graph $G$ to another $(n,m)$-graph $H$ is a vertex mapping that preserves the adjacencies along with their types and directions. The order of 
  a smallest (with respect to the number of vertices) such  $H$ is the $(n,m)$-chromatic number of $G$. 
  Moreover, an $(n,m)$-relative clique $R$ of an $(n,m)$-graph $G$ is a vertex subset of $G$ for which no two distinct vertices of $R$ get identified under 
		any homomorphism of $G$. 
		The $(n,m)$-relative clique number of $G$, denoted by $\omega_{r(n,m)}(G)$, is the maximum $|R|$ such that $R$ is an $(n,m)$-relative clique of $G$. In practice, $(n,m)$-relative cliques are 
  often used for establishing lower bounds of $(n,m)$-chromatic number of graph families.

		Generalizing an open problem posed by Sopena [Discrete Mathematics 2016] in his latest survey on oriented coloring, Chakroborty, Das, Nandi, Roy and Sen [Discrete Applied Mathematics 2022]
		conjectured that 
  $\omega_{r(n,m)}(G) \leq 2 (2n+m)^2 + 2$ for any triangle-free planar $(n,m)$-graph $G$ and that this bound is tight for all $(n,m) \neq (0,1)$.  
		In this article, we positively settle this conjecture by improving the previous upper bound of 
 $\omega_{r(n,m)}(G) \leq 14 (2n+m)^2 + 2$ to 
 $\omega_{r(n,m)}(G) \leq 2 (2n+m)^2 + 2$, and by finding examples of triangle-free planar graphs that achieve this bound. As a consequence of the tightness proof, we also establish a new lower bound of $2 (2n+m)^2 + 2$ for the $(n,m)$-chromatic number of the family of triangle-free planar graphs. 

 \medskip

	\end{abstract}

\section{Introduction}
 The focus of this work is on the resolution of a recent 
 conjecture proposed Chakroborty, Das, Nandi, Roy and 
 Sen~\cite{chakraborty2021clique}  which is in the context of homomorphisms of  
 $(n,m)$-graphs, introduced in 2000 by Ne\v{s}et\v{r}il and 
 Raspaud~\cite{nevsetvril2000colored}.  The conjecture 
 generalizes an open question posed by Sopena~\cite{sopena2016homomorphisms} in his 
 survey on oriented graphs.

\subsection{Colored mixed graphs}
An $(n,m)$-colored mixed graph, or simply, an $(n,m)$-graph $G$ is a graph with $n$ different types of arcs and $m$ different types of edges.
We denote the set of vertices, arcs, and edges of $G$ by $V(G), A(G)$, and $E(G)$ respectively. Also, we denote the underlying graph of $G$ by $und(G)$. 

In this article, we restrict ourselves to studying only those $(n,m)$-graphs whose underlying graphs are simple graphs. We can capture the notions of simple graphs when $(n,m)= (0,1)$~\cite{west2001introduction}, oriented graphs when $(n,m)=(1,0)$~\cite{klostermeyer2004homomorphisms,sopena1997chromatic,sopena2001oriented,sopena2016homomorphisms}, $2$-edge-colored graphs or signed graphs when $(n,m) = (0,2)$~\cite{montejano2010homomorphisms,naserasr2015homomorphisms,naserasr2021homomorphisms,ochem2017homomorphisms,zaslavsky1982signed}, and $m$-edge colored graphs when $(n,m) = (0,m)$~\cite{alon1998homomorphisms}. Therefore, $(n,m)$-graphs help us investigate well-studied families of graphs in a generalized set-up. 
Moreover, the upcoming definitions of graph homomorphisms, chromatic numbers, and clique numbers for $(n,m)$-graphs truly capture the existing concepts of the same for all values of $(n,m)$, including $(0,1)$. 
In this work though, whenever we use the term $(n,m)$-graphs, we mean it for all values of $(n,m) \neq (0,1)$ unless otherwise stated, for convenience.

\subsection{Homomorphisms, chromatic numbers, and clique numbers}
A \textit{homomorphism} of an $(n,m)$-graph $G$ to another $(n,m)$-graph $H$ is a vertex mapping $f: V(G) \to V(H)$ such that  for any arc (resp., edge) $xy$ in $G$, their images induces an arc (resp., edge)  $f(x)f(y)$ of the same type in $H$. If there exists a homomorphism of $G$ to $H$, then we denote it by $G \to H$. 

Using the concept of homomorphisms of such graphs, Ne\v{s}etril and Raspaud~\cite{nevsetvril2000colored} further presented a generalization to the notion of chromatic number by introducing the $(n,m)$-chromatic number of an $(n,m)$-graph as 
$$\chi_{n,m}(G) := \min \{ |V(H)| : G \to H \}. $$ 
For a family $\mathcal{F}$ of simple graphs, the $(n,m)$-chromatic number can be defined as $$\chi_{n,m}(\mathcal{F}) := \max \{\chi_{n,m}(G) : und(G) \in \mathcal{F}\}. $$
This parameter is studied for the family of graphs having bounded acyclic chromatic numbers,
bounded arboricity and acyclic chromatic number, 
bounded maximum degree, 
sparse graphs,
planar graphs and planar graphs with girth restrictions,  
partial $k$-trees and partial $k$-trees with girth restrictions,  
outerplanar graphs and outerplanar graphs with girth restrictions,
etc. across several papers~\cite{fabila2008lower,LNSS2021,montejano2009chromatic,nevsetvril2000colored}. 

Bensmail, Duffy, and Sen~\cite{bensmail2017analogues} contributed to this line of work by introducing and studying generalizations of the concept of clique number. The generalization ramifies into two parameters in the context of $(n,m)$-graphs. 
An \textit{$(n,m)$-relative clique} 
$R \subseteq V(G)$ is a vertex subset satisfying
$|f(R)| = |R|$ for all homomorphisms $f$ of $G$ to any $H$. 
The $(n,m)$-relative clique number of a graph $G$ denoted by $\omega_{r(n,m)}(G)$, is,
$$\omega_{r(n,m)}(G) = \max \{|R|  : R \text{ is an } (n,m)\text{-relative clique of}~ G \}.$$
An \textit{$(n,m)$-absolute clique} 
$A \subseteq G$ is a subgraph of $G$ satisfying
$\chi_{n,m}(A) = |V(A)|$.
The $(n,m)$-absolute clique number of $G$ denoted by $\omega_{a(n,m)}(G)$, is,
$$\omega_{a(n,m)}(G) = \max \{ |V(A)|: A \text{ is an } (n,m)\text{-an absolute clique of}~ G \}.$$  For family $\mathcal{F}$ of simple graphs, both the parameters are defined similarly like $(n,m)$-chromatic number. That is,
$$p(\mathcal{F}) := \max \{p(G) : und(G) \in \mathcal{F}\} $$
where $p \in \{\omega_{r(n,m)}, \omega_{a(n,m)}\}$.  
The two parameters are primarily being studied for 
graphs with bounded degrees, planar, partial $2$-trees and outerplanar graphs, and their subfamilies due to girth restrictions~\cite{bensmail2017analogues,chakraborty2021clique}.

\subsection{Context and motivation}
Observe that the above defined three parameters trivially satisfy the following relation~\cite{bensmail2017analogues}
$$\omega_{a(n,m)}(G) \leq \omega_{r(n,m)}(G) \leq \chi_{n,m}(G).$$
Therefore, the study of one impacts the other. 
The best known lower and upper bounds for the $(n,m)$-chromatic number of planar graphs is given below. The lower bounds for $(n,m)$-chromatic number of planar graphs was given by Fabila-Monroy, Flores, Huemer and Montejano\cite{fabila2008lower} while the upper bound was due to Ne\v{s}et\v{r}il and Raspaud~\cite{nevsetvril2000colored}. Let the family of planar graphs with girth at least $g$ is denoted by $\mathcal{P}_g$.

\begin{theorem}[\cite{fabila2008lower,nevsetvril2000colored}]\label{mixed ch no}
	For the family of planar graphs $\mathcal{P}_3$ we have,
	\begin{align*}
		(2n+m)^3 + 2 (2n+m)^2 + (2n+m) + 1 \leq  \chi_{n,m}(\mathcal{P}_3) \leq 5(2n+m)^4, &\text{ for $m > 0$ even }\\
		(2n+m)^3 +  (2n+m)^2 + (2n+m) + 1 \leq  \chi_{n,m}(\mathcal{P}_3) \leq 5(2n+m)^4, &\text{ otherwise}.
	\end{align*}
\end{theorem}

The above theorem can be treated as an approximate analogue of the Four-Color Theorem in the context of $(n,m)$-graphs. For triangle-free planar graphs, no such dedicated studies have been made, that is, an approximate analogue of the 
Gr\"{o}tzsch\textquotesingle s Theorem is open till date. On the other hand, it turns out, even though finding clique numbers for these two families, that is, the families of planar and triangle-free planar graphs, is trivial for $(n,m)=(0,1)$, and it is quite a challenging problem for the other values of $(n,m)$. There have been dedicated studies to explore these parameters even for $(n,m)=(1,0)$~\cite{das2018study,das2020relative},  and later for  general values~\cite{bensmail2017analogues,chakraborty2021clique}.

We recall the existing bounds for the two clique numbers for the families of planar and triangle-free graphs to place our work in context. For convenience, we present all these results under one theorem. The bounds for the $(n,m)$-absolute clique number of planar graphs $\mathcal{P}_3$ is proved by Bensmail, Duffy and Sen~\cite{bensmail2017analogues}  while the rest are due to  Chakraborty, Das, Nandi, Roy and Sen~\cite{chakraborty2021clique}.

\begin{theorem}\label{th all}
	For the families $\mathcal{P}_3$ of planar and $\mathcal{P}_4$ triangle-free planar graphs, we have
	\begin{enumerate}[(i)]
			\item	 $3 (2n+m)^2 + (2n+m) + 1 \leq \omega_{a(n,m)}(\mathcal{P}_3) \leq 9 (2n+m)^2 + 2(2n+m) + 2$~\cite{bensmail2017analogues},

		\item	 $3 (2n+m)^2 + (2n+m) + 1 \leq \omega_{r(n,m)}(\mathcal{P}_3) \leq 42 (2n+m)^2 - 11$~\cite{chakraborty2021clique},
		
			\item	 $\omega_{a(n,m)}(\mathcal{P}_4) =  (2n+m)^2 + 2$~\cite{chakraborty2021clique},

		\item	 $(2n+m)^2 + 2 \leq \omega_{r(n,m)}(\mathcal{P}_4) \leq  14(2n+m)^2 + 1$~\cite{chakraborty2021clique}.
	\end{enumerate}
\end{theorem}

If one notices, except the $(n,m)$-absolute clique number of triangle-free planar graphs, exact general values are not known for the other parameters from Theorem~\ref{th all}.

%To share a bit of history, let us mention that 
The inception of the analogue to clique number in this set up was due to Klostermeyer and MacGillivray~\cite{klostermeyer2004analogues} where they defined the notion of absolute clique number for $(1,0)$-graphs. In the same work, they proved 
$$\omega_{a(1,0)}(\mathcal{P}_4) \leq 14 \text{ and }
\omega_{a(1,0)}(\mathcal{P}_3) \leq 36,$$  
and conjectured 
$\omega_{a(1,0)}(\mathcal{P}_3) = 15$. 
Nandy, Sen and Sopena~\cite{nandy2016outerplanar} proved the tight bounds 
$$\omega_{a(1,0)}(\mathcal{P}_4) = 6 \text{ and }  \omega_{a(1,0)}(\mathcal{P}_3) = 15$$ 
and thus, positively settled the above mentioned conjecture.

While proving the conjecture, Nandy, Sen and Sopena~\cite{nandy2016outerplanar} introduced the notion of relative cliques for $(1,0)$-graphs which was a useful tool in their proofs. Then 
relative clique number of $(1,0)$-graphs were studied by Das, Prabhu and Sen~\cite{das2018study} in a focused and systematic manner. 
In particular, they proved 
$$10 \leq \omega_{r(1,0)}(\mathcal{P}_4) \leq 14 \text{ and }  15 \leq \omega_{r(1,0)}(\mathcal{P}_3) \leq 32.$$ 
In the meantime, Sopena~\cite{sopena2016homomorphisms} enquired the exact values of $\omega_{r(1,0)}(\mathcal{P}_3)$ and $\omega_{r(1,0)}(\mathcal{P}_4)$ 
in his survey in 2016 (latest to date) on oriented coloring as open problems. 
Das, Nandi and Sen~\cite{das2020relative} solved the second problem by showing 
$\omega_{r(1,0)}(\mathcal{P}_4) = 10$ while the first question remains open still.

Later, Bensmail, Duffy and Sen~\cite{bensmail2017analogues} generalized the notions of relative and absolute clique numbers for all $(n,m)$-graphs and proved Theorem~\ref{th all}(i). 
They~\cite{bensmail2017analogues} also conjectured that lower bound of Theorem~\ref{th all}(i) is tight, generalizing the result $\omega_{a(1,0)}(\mathcal{P}_3) = 15$ due to Nandy, Sen and Sopena~\cite{nandy2016outerplanar}. 

\begin{conjecture}[\cite{bensmail2017analogues}]\label{conj planar absolute}
	For the family  $\mathcal{P}_4$ of triangle-free planar graphs we have
	$$\omega_{a(n,m)}(\mathcal{P}_3) =  3 (2n+m)^2 + (2n+m) + 1.$$
\end{conjecture}

Extending this direction of research, 
Chakraborty, Das, Nandi, Roy and Sen~\cite{chakraborty2021clique} proved a tight bound for the $(n,m)$-absolute clique number of triangle-free planar graphs via 
Theorem~\ref{th all}(iii). They also provided lower and upper bounds for the relative clique numbers of planar and triangle-free planar graphs in Theorem~\ref{th all}(ii) and~(iv). Moreover, in the same work, they made two conjectures about the exact values of 
$\omega_{r(n,m)}(\mathcal{P}_3)$ and $\omega_{r(n,m)}(\mathcal{P}_4)$, generalizing the open problems proposed by Sopena~\cite{sopena2016homomorphisms} in the context of $(1,0)$-graphs.   

\begin{conjecture}[\cite{chakraborty2021clique}]\label{conj planar relative}
	For the family  $\mathcal{P}_3$ of planar graphs we have
	$$\omega_{r(n,m)}(\mathcal{P}_3) =  3 (2n+m)^2 + (2n+m) + 1.$$
\end{conjecture}

\begin{conjecture}[\cite{chakraborty2021clique}]\label{conj triangle-free planar relative}
	For the family  $\mathcal{P}_4$ of triangle-free planar graphs we have
	$$\omega_{r(n,m)}(\mathcal{P}_4) =  2(2n+m)^2 + 2.$$
\end{conjecture}

\subsection{Our contributions}
 As mentioned before, Das, Nandi and Sen~\cite{das2020relative} solved Conjecture~\ref{conj triangle-free planar relative} for $(n,m)=(1,0)$. Moreover, they noted that 
 the solution for $(n,m) = (0,2)$ is similar. Thus, we may consider the conjecture to be solved when $(2n+m)= 2$. In this article, we  positively settle this conjecture 
when $(2n+m) \geq 3$ by proving the following. 

\begin{theorem}\label{thm main}
	For the family  $\mathcal{P}_4$ of triangle-free planar graphs we have
	$$\omega_{r(n,m)}(\mathcal{P}_4) =  2(2n+m)^2 + 2$$
	for  $(2n+m) \geq 3$. 
\end{theorem}

As a direct consequence, we find a new lower bound for the $(n,m)$-chromatic number of triangle-free planar graphs.

\begin{theorem}\label{thm chi}
	For the family  $\mathcal{P}_4$ of triangle-free planar graphs we have
	$$\chi_{n,m}(\mathcal{P}_4) \geq  2(2n+m)^2 + 2$$
	for  $(2n+m) \geq 3$. 
\end{theorem}

\subsection{The significance of relative cliques}
The focus parameter of this article, $(n,m)$-relative cliques, is not such a popularly studied parameter. 
However, a closer look may imply otherwise, hence the need of this subsection. As $(n,m)$-absolute cliques are special cases of $(n,m)$-relative cliques, we will try to point out works done on both areas to highlight the significance of 
$(n,m)$-relative cliques. Every point written below is to the best of our knowledge. 

\begin{itemize}
    \item For $(n,m)$-graphs, there have been two works dedicated on $(n,m)$-absolute and $(n,m)$-relative cliques: \cite{bensmail2017analogues} and 
    \cite{chakraborty2021clique}. 

    \item For oriented graphs, that is, for $(n,m) = (1,0)$, a number of papers dealt with the notions of absolute and relative cliques: \cite{bensmail2016complexity, das2018study, das2020relative, dybizbanski2020oriented, klostermeyer2004analogues, nandy2016outerplanar}.
    A similar variant of homomorphisms of oriented graphs, called pushable homomorphisms of oriented graphs, also studies the notions of absolute and relative cliques~\cite{bensmail2017oriented}. 

    \item For $2$-edge-colored graphs, that is, for $(n,m) = (0,2)$, a some of papers dealt with the notions of absolute and relative cliques. However, it is usually for a particular variant of homomorphisms of $2$-edge colored graphs, called homomorphisms of signed graphs, which recently gained popularity due to its underlying connections with the theory of graph coloring. The notions of absolute and relative cliques exists, and are relevant,  for signed graphs also. The following are some works on them: \cite{bensmail2020classification, 
    das2019relative, das2020relative-signed,
    naserasr2015homomorphisms}.    
\end{itemize}

The $(n,m)$-relative (resp., absolute) clique numbers are natural lower bounds for $(n,m)$-chromatic numbers. In practice, one of the most popularly used 
    methods to establish a lower bound for 
    $(n,m)$-chromatic number of a family of graphs is to  provide an example of a graph from that family having high $(n,m)$-relative clique number. This trick has been used even before the notions of absolute and relative cliques came into existence. This is possibly the most important significance of the notions of absolute and relative cliques.  In the following we list such instances.

\begin{itemize}
    \item For $(n,m)$-graphs, lower bounds for 
    $(n,m)$-chromatic number of several graph families were established using 
    examples having high $(n,m)$-relative clique number across the following works: \cite{fabila2008lower, 
    montejano2009chromatic, nevsetvril2000colored}. 

    \item For oriented graphs, that is, for $(n,m) = (1,0)$, lower bounds for 
    $(1,0)$-chromatic number (popularly known as oriented chromatic number) of several graph families were established using 
    examples having high $(n,m)$-relative clique number across a number of works. Most of them can be 
    found in the latest survey by Sopena~\cite{sopena2016homomorphisms}. We would like to make a special mention to~\cite{marshall2007homomorphism, marshall2015oriented}, where relative cliques and the analysis of their interactions was used to improve the lower bound for the oriented chromatic number of the family of planar graphs. 
    Similar phenomenon can be observed for pushable chromatic number of oriented graphs in the following works: \cite{bensmail2021pushable, klostermeyer2004homomorphisms}. 

    \item For $2$-edge colored graphs, that is, for $(n,m) = (0,2)$,  
    $(0,2)$-chromatic number\footnote{This parameter has been studied under several different names across the literature, such as: $2$-edge-colored chromatic number, signified chromatic number, and signed-preserving chromatic number.} 
    and of chromatic number of signed graphs, the use of absolute and relative cliques to establish lower bounds can be observed  in the following works: \cite{bensmail2022signed, naserasr2015homomorphisms, ochem2017homomorphisms}. 
\end{itemize}

In particular for oriented graphs, that is, for $(n,m) = (1,0)$, which is also historically the oldest and most studied specific case of homomorphisms of $(n,m)$-graphs, absolute and relative cliques appear in two other significant areas. 

\begin{itemize}
    \item Borodin, Fon-Der-Flaass, Kostochka, Raspaud, and Sopena~\cite{borodin2001deeply} introduced and studied the notion of deeply critical oriented graphs and asked the existence of ``deeply critical $(1,0)$-absolute cliques'' (even before the term absolute clique came into existence) on prescribed number,  of vertices. They themselves proved such existence when the number of vertices are of the form $2 \cdot 3^n -1$. Later, answering to an open question asked in the above mentioned paper, Duffy, Pavan, Sandeep, and Sen~\cite{pavan} found examples 
 of such $(1,0)$-absolute cliques for all odd values of number of vertices, except $1$, $3$, and $7$.

\item In relation to the famous degree-diameter problem, Erd\H{o}s, R\'{e}nyi and S\'{o}s~\cite{erdos1966problem} defined the function $$h_d(n)=\min\{|E(G)| : diam(G) \leq d \text{ and }|V(G)|=n\}$$ as the minimum number of edges among graphs having diameter at most $d$. Its directed (oriented) analogue turns out to be very interesting even 
for $d=2$. In fact, the oriented graphs having (weak) diameter $2$ are exactly the $(1,0)$-absolute cliques. A number of works are dedicated in quest of finding the exact description of this function~\cite{furedi1998minimal,katona1967problem,kostochka1999minimum}. 
\end{itemize}

Thus, $(n,m)$-relative cliques are important objects of study in the context of homomorphisms of $(n,m)$-graphs.

\subsection{Organization of the article}
This article is structured as follows. In Section~\ref{sec prelims} we give the necessary definitions and notations used in this article. In Section~\ref{sec proof} we prove  Theorem~\ref{thm main}. In Section~\ref{sec conclusions} 
we share our concluding remarks.

\medskip

\noindent \textbf{Note:} A preliminary version of this article has been published in IWOCA 2022~\cite{nandi2022relative}.

\section{Preliminaries}\label{sec prelims}
In this section, we introduce some notation to help us write the proofs.  Also, we follow West~\cite{west2001introduction} for standard definitions, notation and terminology. 

Given an $(n,m)$-graph $G$ the different types of arcs in $G$ are distinguished by $n$ different labels $2, 4, \cdots, 2n$. 
To be precise, an arc $xy$ with label $2i$ is an arc 
of type $i$ from $x$ to $y$. 
In fact, in such a situation, $y$ is called a $2i$-neighbor of $x$, or equivalently, $x$ is called a $(2i-1)$-neighbor of $y$. The different types of edges in $G$ are distinguished 
by $m$ different types of labels $2n+1, 2n+2, \cdots, 2n+m$. 
Here, an edge $xy$ with label $2n+j$ is an edge of type $2n+j$ between $x$ and $y$. In this case, $x$ and $y$ are called $(2n+j)$-neighbors of each other.

Usually, throughout the article, we use the Greek alphabets such as $\alpha, \beta, \gamma$ and their variants (like $\alpha'$, $\beta'$, $\gamma'$, etc.) to denote these labels. Therefore, whenever we use such a symbol, say $\alpha$, one may assume it to be an integer between $1$ and $(2n+m)$. Let 
us denote the set of all $\alpha$-neighbors of $x$ as $N^{\alpha}(x)$. Two vertices $x, y$ \textit{agree} on a third vertex $z$ if $z \in N^{\alpha}(x) \cap N^{\alpha}(y)$ for some $\alpha$, and \textit{disagree} on $z$ otherwise. 

A \textit{special $2$-path} is a $2$-path $uwv$ such that $u,v$ disagrees with each other on $w$. 
In an $(n,m)$-graph $G$, a vertex $u$ \textit{sees} a vertex $v$ if they are either adjacent, or are connected by a special $2$-path. If $u$ and $v$ are connected by a special $2$-path with $w$ as the \textit{internal vertex}, then it is said that $u$ sees $v$ via $w$ or equivalently, 
$v$ sees $u$ via $w$.
We recall a useful characterization of an $(n,m)$-relative clique due to 
Bensmail, Duffy and Sen~\cite{bensmail2017analogues}.   

\begin{lemma}[\cite{bensmail2017analogues}]
	Two distinct vertices of an $(n,m)$-graph $G$ are part of a relative clique if and only if they are either adjacent or 
	connected by a special 2-path in $G$, that is, they see each other.
\end{lemma}

\section{Proof of Theorem~\ref{thm main}}\label{sec proof}
The proof is lengthy and involved. So for the reader's convenience, we have presented it under separate indicative headers and have distributed it into different propositions, lemmas, and observations.

\subsection{Lower bound}
The lower bound of Theorem \ref{thm main} is given by the following proposition.

\begin{proposition}
	There exists a triangle-free planar  $(n,m)$-graph $G$ such that
	$$\omega_{r(n,m)}(G) = 2(2n+m)^2 + 2.$$ 
\end{proposition}

\begin{proof}
	We construct an example to prove this claim. We start with a $K_{2,2(2n+m)^2}$. Let the vertices in the partite set consisting of 
	only two vertices be $x$ and $y$. Let the vertices from the other partite set be 
	$$\{g_{\alpha \beta}, g'_{\alpha \beta}: \text{ for all } \alpha, \beta \in \{1,2,  \cdots , 2n+m \} \}.$$  
	Next, assign colors and directions to the  edges in such a way that
	$g_{\alpha \beta}$ and $g'_{\alpha \beta}$ become 
	$\alpha$-neighbors of $x$ and $\beta$-neighbors of $y$. 
	Finally, add special $2$-paths 
	of the form $g_{\alpha \beta}h_{\alpha \beta} g'_{\alpha \beta}$ connecting $g_{\alpha \beta}$
	and $g'_{\alpha \beta}$. 
	
	Notice that, $g_{\alpha \beta}$ (or $g'_{\alpha \beta}$) sees
	$g_{\alpha' \beta'}$ (or $g'_{\alpha' \beta'}$) via $x$ if $\alpha \neq \alpha'$ or via $y$ if $\beta \neq \beta'$.
	If $\alpha = \alpha'$ and $\beta = \beta'$, then they are the two vertices
	$g_{\alpha \beta}$ and $g'_{\alpha \beta}$ and they see each other via 
	$h_{\alpha \beta}$. 
	On the other hand, $x$ and $y$ are adjacent to each $g_{\alpha \beta}$, while they see each other via $g_{\alpha' \beta'}$ for some $\alpha' \neq \beta'$. Thus, the $g_{\alpha \beta}$\textquotesingle s together with $x$ and 
	$y$, form an $(n,m)$-relative clique of cardinality $2(2n+m)^2 + 2$. 
\end{proof}

\subsection{Strategy and set up for the upper bound}
Next, we concentrate on the upper bound. We prove the upper bound by the method of contradiction. Therefore, let us assume that $\omega_{r(n,m)}(\mathcal{P}_4) > 2(2n+m)^2+2$. 

Let $n_3(G)$ be the number of vertices with a degree at least $3$. 
A \textit{critical $(n,m)$-relative clique} $H$ 
for the family $\mathcal{P}_4$ of triangle-free planar graphs 
is an $(n,m)$-graph $H$ satisfying the following properties:

\begin{enumerate}[(i)]
	\item $und(H) \in \mathcal{P}_4$,

	\item $\omega_{r(n,m)}(H) = \omega_{r(n,m)}(\mathcal{P}_4),$

	\item $\omega_{r(n,m)}(H^*) < \omega_{r(n,m)}(H)$, if $$(n_3(H^*), |V(H^*)|,|E(und(H^*))|) < (n_3(H), |V(H)|,|E(und(H))|)$$
	in the dictionary ordering, where $H^* \in \mathcal{P}_4.$
\end{enumerate}

We consider a particular critical $(n,m)$-relative clique $H$ for the rest of this proof. As $H$ is a triangle-free planar graph, we fix a particular planar embedding of it. Whenever we deal with something dependent on the embedding of $H$ or a portion of it, we refer to this particular fixed embedding. We also fix a particular $(n,m)$-relative clique $R$ of cardinality $\omega_{r(n,m)}(H)$. Furthermore, the vertices of $R$ are called \textit{good vertices} and those of $S = V(H) \setminus R$ are called the \textit{helper vertices}. Notice that, 
$$\omega_{r(n,m)}(H) = |R| = \omega_{r(n,m)}(\mathcal{P}_4) > 2(2n+m)^2+2$$
due to our basic assumption. Thus, the goal of the proof is  to contradict this strict inequality.

\subsection{Key lemmas}
We recall some useful results due to 
Chakraborty, Das, Nandi, Roy and Sen~\cite{chakraborty2021clique}.

\begin{lemma}[\cite{chakraborty2021clique}] \label{lem Sind}
	The $(n,m)$-graph $H$ is connected and the set $S$ of helper vertices in $H$ is an independent set. 
\end{lemma}

We recall a modified version of another lemma due to Chakraborty, Das, Nandi, Roy and Sen~\cite{chakraborty2021clique} which we use for our proofs. 

\begin{lemma}[\cite{chakraborty2021clique}]\label{lem 2p alpha max}
	A vertex $x \in V(H)$ can have at most $2(2n+m)$ good $\alpha$-neighbors in $H$ for any $\alpha$. 
\end{lemma}

\begin{figure}
	\begin{center}
		\begin{tikzpicture}  
			[scale=0.9,auto=center]

			\draw (0,8) --node[above]{$\alpha$} (-3,6);
		%	\draw (0,6) --node[pos = 0.35, below right]{$\alpha$} (0,4);
		\draw (0,8) -- node[above right]{$\alpha$} (0,6); 
			\draw (0,8) --node[above]{$\alpha$} (3,6);
			\draw (0,4) --node[below]{$\beta$} (-3,6);
			\draw (0,4) --node[below right]{$\beta$} (0,6);
			\draw (0,4) --node[below]{$\beta$} (3,6);
			\draw[] (-3,6) -- (-1.5,6) -- (0,6) -- (1.5,6) -- (3,6);
			
			\draw[fill] (-1.5,6) circle[radius=0.08];
			 \draw [ >=stealth] (0,2.5) to [bend left=45] (-3,6) ;
			  \draw [ >=stealth] (3,6) to [bend left=45] (0,2.5) ;

			\node at (-1.5,5.6) {$h_1$};
				\node at (1.5,5.6) {$h_2$};
					\node at (0,2.2) {$h_3$};
			\node[circle,draw,fill=black!20] at (0,8) {$x$};
			\node[circle,draw,fill=black!20] at (-3,6) {$g_1$};
				\draw[fill] (-1.5,6) circle[radius=0.08] ;
		\draw[fill] (1.5,6) circle[radius=0.08];
			\draw[fill] (0,2.5) circle[radius=0.08] ;
		
			\node[circle,draw,fill=black!20] at (0,6) {$g_2$};
		
			\node[circle,draw,fill=black!20] at (3,6) {$g_3$};
			
			\node[circle,draw,fill=black!20] at (0,4) {$y$};

		\end{tikzpicture}  
		
	\end{center}
	
	\caption{The exceptional configuration of Lemma~\ref{lem 2 good in alpha-beta}. The symbols $\alpha$ and $\beta$ are adjacency types with respect to $x, y$, respectively and $h_1, h_2,h_3$ denote the helper vertices.}
	\label{fig good vertices}
\end{figure}
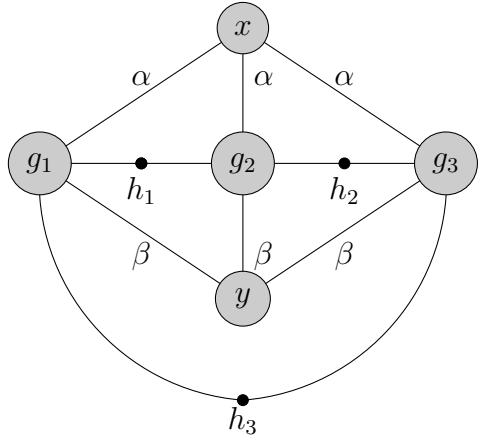
We prove a similar result for common good neighbors of two distinct vertices  as well, but for an exception. Notice that, it is possible for two vertices $x,y$ of $H$ to agree on three good vertices as depicted in Fig.~\ref{fig good vertices} where the good vertices are also able to see each other via distinct helpers. 
Turns out this is the only exception, otherwise, 
two distinct vertices of $H$ can agree on at most $2$ good vertices.

\begin{lemma}\label{lem 2 good in alpha-beta}
	Two distinct vertices $x,y$ of $H$ can agree on at most $2$ good vertices, except when $x,y$ agree on exactly three good vertices which are mutually connected by internally vertex disjoint special two paths as in Fig.~\ref{fig good vertices}. 
\end{lemma}

\begin{proof}
	Suppose $x$ and $y$ have three good neighbors $g_1, g_2, g_3$ from $N^{\alpha}(x) \cap N^{\beta}(y)$, and at least another common neighbor $g_4$. Also assume that $g_1, g_2, g_3, g_4$ are arranged in an 
	anti-clockwise order around $x$ in the fixed planar embedding of $H$. 
	Notice that, as $H$ is triangle-free, $g_1$ must see $g_3$ via some 
	$h \not\in \{x, y, g_1, g_2, g_3, g_4\}$. This is not possible to achieve keeping the graph planar. Hence, if $x$ and $y$ have four common good neighbors, then it is not possible for three of their neighbors to agree 
	on $x$ and $y$. 
	
	If $x$ and $y$ have exactly three common neighbors $g_1, g_2, g_3$, and all three agree on $x$ and $y$, then as $H$ is triangle-free, $g_i$ must see $g_j$ via some $h_{ij}$ for all $i < j$.  These $h_{ij}$\textquotesingle s must be distinct, otherwise a $K_{3,3}$ will be created. Hence Fig.~\ref{fig good vertices} is forced.
\end{proof}

In general, two vertices in $H$ can have at most $2(2n+m)^2$ many common good neighbors.

\begin{lemma}\label{lem common neighbors upper bound}
	Two distinct vertices $x,y$ of $H$ can have at most $2(2n+m)^2$ good vertices in their common neighborhood.  
\end{lemma}

\begin{proof}
	As $(2n+m) \geq 2$, and thus, $2(2n+m)^2 \geq 8$, the exceptional case mentioned in Lemma~\ref{lem 2 good in alpha-beta} satisfies the condition of the statement of this lemma. 
	
	In other case, Lemma~\ref{lem 2 good in alpha-beta} implies that 
	$$|N^{\alpha}(x) \cap N^{\beta}(y) \cap R| \leq 2.$$
	As $\alpha, \beta$ can be chosen from a set of $(2n+m)$ integers, we can have a total of $(2n+m)^2$ distinct $(\alpha, \beta)$ pairs. Therefore, there can be at most $2(2n+m)^2$ good vertices in 
	$N(x) \cap N(y)$. 
\end{proof}

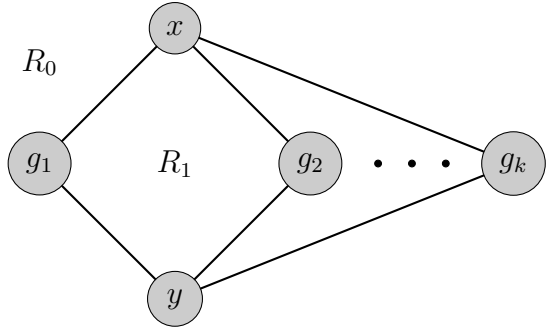
\begin{figure}
	\begin{center}
		\begin{tikzpicture}  
			[scale=.9,auto=center] % here, node/.style is the style pre-defined, that will be the default layout of all the nodes. You can also create different forms for different nodes.  

			\draw[thick] (0,4) -- (-2,2);
			\draw[thick] (0,4) -- (2,2);
			\draw[thick] (0,4) -- (5,2);
			\draw[thick] (0,0) -- (-2,2);
			\draw[thick] (0,0) -- (2,2);
			\draw[thick] (0,0) -- (5,2);

			\node[circle,draw,fill=black!20] at (0,4) {$x$};
			\node[circle,draw,fill=black!20] at (0,0) {$y$};
			\node[circle,draw,fill=black!20] at (-2,2) {$g_1$};
			\node[circle,draw,fill=black!20] at (2,2) {$g_2$};
			\node[circle,draw,fill=black!20] at (5,2) {$g_k$};
			\node at (-2,3.5) {$R_0$};
			\node at (0,2) {$R_1$};
			\draw[fill] (3,2) circle[radius=0.05];
			\draw[fill] (3.5,2) circle[radius=0.05];
			\draw[fill] (4,2) circle[radius=0.05];

		\end{tikzpicture}  
		
	\end{center}
	\caption{The configuration $\mathbb{F}_k$. }
	\label{fig Fk}
\end{figure}

\subsection{The forbidden configuration $\mathbb{F}_k$}
The configuration $\mathbb{F}_k$ consists of two vertices $x,y$ of $H$ which may not necessarily be good vertices and $k$ common good neighbors of $x,y$, namely, $g_1, g_2, \cdots, g_k$. Also, we assume a default embedding
of $\mathbb{F}_k$ where
$g_1, g_2, \cdots, g_k$ are arranged in an anti-clockwise order around $x$.  
Without loss of generality, 
we may assume the embedding of $\mathbb{F}_k$, induced by the fixed embedding of $H$, to be such that the boundary of the unbounded region is the $4$-cycle $xg_1yg_kx$. Thus the $4$-cycle $xg_1yg_kx$ divides the plane into two regions: a bounded region containing other $g_i$\textquotesingle s called $R$, and an unbounded region called $R_0$. Moreover, the 
$2$-paths of the type $xg_iy$ further subdivides the region $R$ into 
$(k-1)$ regions. These regions are called $R_1, R_2, \cdots, R_{k-1}$, where 
$R_i$ is bounded by $xg_iyg_{i+1}x$, for $i \in \{1,2,  \cdots, k-1\}$. 
See Fig.~\ref{fig Fk} for a pictorial representation of this planar embedding.

If it is not possible for $H$ to contain a $\mathbb{F}_k$ for some $k$, then we say that $\mathbb{F}_k$ is \textit{forbidden} in $H$. Therefore, Lemma~\ref{lem common neighbors upper bound} essentially says that $\mathbb{F}_k$, for all $k \geq 2(2n+m)^2 +1$ is forbidden. Next, we show that $\mathbb{F}_k$, for all $k \geq 3$ is forbidden as a key step of our proof.

\begin{lemma}\label{lem Fk}
	The configuration $\mathbb{F}_k$, for all $k \geq 3$ is forbidden in $H$. 
\end{lemma}

\begin{proof}
	We prove this by strong backward induction. Suppose that the statement is true for all $k \geq t+1$ and we prove that it is  true for $k = t$, where $t \geq 3$. The base case is taken care of by Lemma~\ref{lem common neighbors upper bound}, where it is proved that $\mathbb{F}_k$ is forbidden for all $k \geq 2(2n+m)^2 +1$ in $H$.

	Now, notice that, if a good vertex $g$ which is not part of $\mathbb{F}_k$, belongs to region $R_0$ (resp., $R_1$, $R_2$), then it must see $g_2$ (resp., $g_3$, $g_1$)  via $x$ or $y$. 
	Observe that, it is not possible for $g$ to be adjacent to both $x$ and $y$, as otherwise, it will create a 
	$\mathbb{F}_{t+1}$, which is forbidden due to the induction hypothesis.

	Therefore, every good vertex of $H$, other than $x,y$, is either adjacent to both $x$ and $y$ or exactly one of them. If in case, all of them are adjacent to $x$, then we can say that there are maximum 
	$2(2n+m)^2$ good vertices in the neighborhood of $x$ due to Lemma~\ref{lem 2p alpha max}. 
	This implies that $H$ has at most $2(2n+m)^2+2$ good vertices, the additional two vertices coming from counting $x$ and $y$.  
	Thus it is not possible for all good neighbors other than $x,y$ to be adjacent to $x$. Similarly, it is not possible for all good neighbors other than $x,y$ to be adjacent to $y$.

	Therefore, there must be at least one good vertex which is adjacent to $x$ but not to $y$ and at least one good vertex that is adjacent to $y$ but not to  $x$. 
	Notice that, these vertices must belong to the same region $R_i$ for some $i \in \{0,1, \cdots, R_{k-1}\}$, as otherwise they will not be able to see each other. Hence, without loss of generality, we may assume that all good vertices other than the ones contained in 
	$\mathbb{F}_k$ belong to $R_1$. 
	These good vertices are adjacent to exactly one of $x$ and $y$. 
	Let us assume that $x_1, x_2, \cdots , x_p$ are the good vertices adjacent to $x$ and not to $y$ 
	while $y_1, y_2, \cdots, y_q$ are the good vertices adjacent to $y$ and not to $x$. 
	Also suppose that 
	$x_1, x_2, \cdots , x_p$ are arranged in an anti-clockwise order around $x$ and 
	$y_1, y_2, \cdots, y_q$ are arranged in a clockwise order around $y$. 
	Also for convenience, 
	$x_i$\textquotesingle s are called private good neighbors of $x$ and $y_j$\textquotesingle s are called private good neighbors of $y$. 
	Similarly, $g_i$\textquotesingle s are called common good neighbors of $x$ and $y$. See Fig.~\ref{fig private neighbors} for a pictorial reference. 
	Furthermore, due to symmetry, we may assume that $p \geq q \geq 1$.

	\begin{figure}[t]
		\begin{center}
			\begin{tikzpicture}  
				[scale=.5,auto=center] % here, node/.style is the style pre-defined, that will be the default layout of all the nodes. You can also create different forms for different nodes.  
				
				\node at (-1.5,3.7) {$x_1$};
				\node at (-0.4,3.7) {$x_2$};
				\node at (1.8,3.7) {$x_p$};
				\node at (-1.5,0.3) {$y_1$};
				\node at (-0.4,0.3) {$y_2$};
				\node at (1.6,0.3) {$y_q$};
				
				\draw[fill] (-1.5,4) circle[radius=0.05];
				\draw[fill] (-0.4,4) circle[radius=0.05];
				\draw[fill] (1.5,4) circle[radius=0.05];   
				\draw[fill] (0,4) circle[radius=0.05];
				\draw[fill] (0.5,4) circle[radius=0.05];
				\draw[fill] (1,4) circle[radius=0.05]; 
				\draw[fill] (-1.5,0) circle[radius=0.05];
				\draw[fill] (-0.4,0) circle[radius=0.05];
				\draw[fill] (1.6,0) circle[radius=0.05]; 
				\draw[fill] (0.1,0) circle[radius=0.05];
				\draw[fill] (0.5,0) circle[radius=0.05];
				\draw[fill] (1,0) circle[radius=0.05]; 
				
				\draw[thick] (0,6) -- (-4,2);
				\draw[thick] (0,6) -- (4,2);
				\draw[thick] (0,6) -- (8,2);
				\draw[thick] (0,6) -- (-1.5,4);
				\draw[thick] (0,6) -- (-0.4,4);
				\draw[thick] (0,6) -- (1.5,4);
				\draw[thick] (0,-2) -- (-4,2);
				\draw[thick] (0,-2) -- (4,2);
				\draw[thick] (0,-2) -- (8,2);
				\draw[thick] (0,-2) -- (-1.5,0);
				\draw[thick] (0,-2) -- (-0.4,0);
				\draw[thick] (0,-2) -- (1.6,0);

				\node[circle,draw,fill=black!20] at (0,6) {$x$};
				\node[circle,draw,fill=black!20] at (0,-2) {$y$};
				\node[circle,draw,fill=black!20] at (-4,2) {$g_1$};
				\node[circle,draw,fill=black!20] at (4,2) {$g_2$};
				\node[circle,draw,fill=black!20] at (8,2) {$g_k$};

				\draw[fill] (5,2) circle[radius=0.05];
				\draw[fill] (6,2) circle[radius=0.05];
				\draw[fill] (7,2) circle[radius=0.05];

			\end{tikzpicture}  
			
		\end{center}
		\caption{The private good neighbors of $x$ and $y$. }
		\label{fig private neighbors}
	\end{figure}
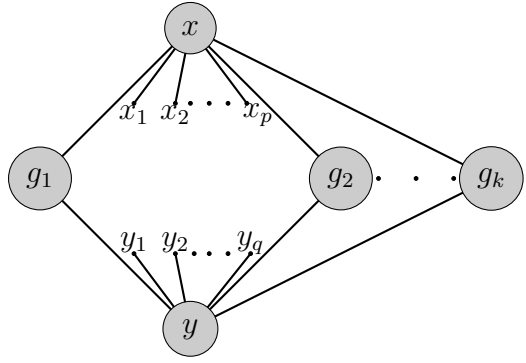

\medskip

	Observe that, if $x$ has a  private good $\alpha$-neighbor, then it cannot have a   
	common good $\alpha$-neighbor other than, possibly, $g_1, g_2$. 
	A similar statement holds for $y$ and its private neighbors as well. 
	Thus, 
	if $x$ (resp., $y$) has $a$ (resp., $b$) different types of 
	adjacencies with its private neighbors, 
	then $x, y$ can have at most $(2n+m - a)(2n+m - b)$ types (with respect to adjacencies with $x$ and $y$) of common neighbors (other than, possibly,  $g_1, g_2$). Thus, due to Lemma~\ref{lem 2 good in alpha-beta}, assuming $\mathbb{F}_k$ is not the exceptional case, there are at most 
	$2(2n+m - a)(2n+m - b)$ common good neighbors of $x,y$ (other than, possibly, $g_1, g_2$). 
	On the other hand, due to Lemma~\ref{lem 2p alpha max}, there are at most $2(2n+m)a$ private good neighbors of 
	$x$ and $2(2n+m)b$ private good neighbors of $y$. Notice that, it is possible to count $g_1$ (resp., $g_2$) 
	among the number of common good neighbors of $x,y$ if it is neither adjacent to $x$ with any of the 
	$a$ types of adjacencies, nor to $y$ with any of the $b$ types of adjacencies. Otherwise, it is possible to count them along with the private neighbors.  Thus the total number of good neighbors in $H$ is at most
	$$ 2(2n+m - a)(2n+m - b) +  2(2n+m)a +2(2n+m)b +2 =  2(2n+m)^2 + 2ab + 2.$$
Even though the above equation does not provide the bound we want, we will consider three different cases and in each case a modified version of the above equation will provide us the required bound.

\medskip

	The above equation gets modified when $p \geq q \geq 2$. In that case note that if $x_1$ is adjacent to $y_q$, then it is not possible for $x_p$ to see $y_1$. Hence, $x_1$ must see $y_q$ via some 
	$h$. This will force every $x_i$ to see $y_j$ via $h$, except for, maybe when $i=j=1$ or when $i=p, j=q$. 
	If $x_1$ agrees with $y_1$ on $h$, then it is not possible for any other $x_i$ (resp., $y_j$), 
	for $i, j \neq 1$,  to agree with 
	$x_1$ on $h$, as otherwise, they cannot see $y_1$ (resp., $x_1$). 
	A similar statement holds if $x_p$ agrees with $y_q$ on $h$. 
	Suppose if there are $c$ types of adjacencies of $h$ among private neighbors of $x$ and 
	$d$ types of adjacencies among private neighbors of $y$, the value of $p+q$ can be at most 
	$2ac +2bd$.
	Therefore, without loss of generality assuming $a \geq b$, 
	the revised upper bound of $p+q$ is 
	$$p+q \leq 2ac +2bd \leq 2a(c+d) \leq 2(2n+m)a.$$
	Hence, the revised upper bound for the total number of 
	good neighbors in $H$ is 
	$$ 2(2n+m - a)(2n+m - b) +  2(2n+m)a + 2 \leq 2(2n+m)^2 +2.$$
	Thus, we have a contradiction in the case when $p \geq q \geq 2$.

 \medskip 
	
	Next, let us concentrate on the case when 
	$q=1$. In this case, each $x_i$ must see $y$ via $y_1$ or via some $h_i$. However, if $x_i$ sees $y$ via
	some $h_i$ and is non-adjacent to $y_1$, then it must see $y_1$ via some $h'_i$. 
	Hence, in the worst case scenario, $x_1$ and $x_p$ will 
	see $y$ and $y_1$ via $h_1, h_p$ and 
	$h'_1, h'_p$, respectively while the other $x_i$'s will see $y$ via $y_1$. 
	In this case, unless $p \leq 2$, 
	it is not possible for $x_1$ to agree with $x_p$ 
	on $x$. 
	Also, $x_i$'s (for $i \neq 1, p$) disagree with $y$ on $y_1$. Therefore,  $p+q \leq 2(2n+m)a + 1$ in this scenario. Hence, when $q=1$ and $p \geq 3$, the total number of good vertices in $H$ is at most
      $$ 2(2n+m - a)(2n+m - 1) +  2(2n+m)a + 1 + 2 \leq 2(2n+m)^2 + 3 - 2 (2n+m - a). $$
Observe that $a < 2n+m$, as none of the private neighbors of $x$ can agree with $g_3$ on $x$. Hence,  $(2n+m - a) \geq 1$, which implies that the total number of good vertices in $H$ is at most $2 (2n+m)^2 + 2.$
	
Thus, we have a contradiction in the case when $p \geq 3$ and $q = 1$. 	

 \medskip 
 
 On the other hand, if $p \leq 2$, then the total number of private neighbors is at most $3$. 
	In this scenario,  the total number of good vertices in $H$ is at most 
	$$ 2(2n+m - a)(2n+m - 1) +  3 + 2 \leq 2(2n+m)^2 - 2(2n+m)(a+1) + a + 3 +2.$$
	Notice that the above bound is at most $2(2n+m)^2 + 2$ as $2n+m \geq 2$ and $a = 1$ or $2$. 
Thus, we have a contradiction in the case when $p \leq 2$.
\medskip  	
 
 This completes the proof, except when $\mathbb{F}_k$ is the exceptional case.  
If $\mathbb{F}_k$ is  the exceptional case, that is, the one described in Fig.~\ref{fig good vertices}, and assume that the common neighbors are $\alpha$-neighbors of $x$ and $\beta$-neighbors of $y$.   
	Then $x$ (resp., $y$) cannot have any private $\alpha$-neighbor (resp. $\beta$-neighbor) as they cannot see all the three common neighbors. 
 Therefore, $a,b \leq 2n+m - 1$. As there are exactly 
 three common neighbors, the total number of good vertices 
in $H$ is at most
$$ 2(2n+m - 1)(2n+m - 1) +  3 + 2 \leq 2(2n+m)^2 -4(2n+m) +7.$$
 As $(2n+m) \geq 3$, the above bound is at most $2(2n+m)^2 + 2$. This is a contradiction, and this completes the proof. 
\end{proof}

\subsection{Proving the existence of a good vertex with low degree}
We now restrict the number of good neighbors for any vertex of $H$. To do so, we need some supporting results.

%\begin{lemma}\label{lem at most four alpha}
%It is not possible for a vertex $x \in V(H)$ to have four or more good 
%$\alpha$-neighbors  for $\alpha \in \{-n, -(n-1), \cdots, n+m\}$. 
%\end{lemma}
%
%\begin{proof}
%Suppose $x$ has at least four good $\alpha$-neighbors $g_1, g_2, g_3, g_4$ arranged in an anti-clockwise order around $x$ in the fixed embedding of $H$. Observe that, as $H$ is triangle-free graph, the only way for $g_1$ to see $g_3$ is via some vertex $h \not\in \{x, g_1, g_2, g_3, g_4\}$. Notice that the 
%$4$-cycle $xg_1hg_3x$ divides the plane into two connected regions: the bounded one containing $g_2$ and the unbounded one containing $g_4$. Notice that, the only way for $g_2$ to see $g_4$ is via $h$. 
%This will result in a $\mathbb{F}_4$, which is a forbidden configuration for $H$ 
%due to Lemma~\ref{lem Fk}.
%\end{proof}

\begin{lemma}\label{lem at most 3 alpha+2beta}
	It is not possible for a vertex $x \in V(H)$ to have four or more good 
	neighbors, with at least three of them being $\alpha$-neighbors.
\end{lemma}

\begin{proof}
	Suppose $x$ has at least four good neighbors $g_1, g_2, g_3, g_4$ 
	arranged in an anti-clockwise order around $x$ in 
	the fixed embedding of $H$, where $g_1, g_2, g_3$ are $\alpha$-neighbors.
	
	Observe that, as $H$ is a triangle-free graph, 
	the only way for $g_1$ to see $g_3$ is via 
	some vertex $h \not\in \{x, g_1, g_2, g_3, g_4\}$.
	Similarly, $g_2$ must see $g_i$ via some $h_i$, for all $i \in \{1, 3\}$. 
	Notice that, $h, h_1, h_3$ are all distinct vertices, as otherwise, a 
	$\mathbb{F}_3$ will be created 
	which is forbidden due to Lemma~\ref{lem Fk}. Also
	as there is no way for $h_1$ and $h_3$ to see $g_4$, thus 
	$h_1, h_3$ cannot be good vertices.

	On the other hand, the $4$-cycle $xg_1hg_3x$ divides the plane into two connected regions: the bounded region containing $g_2$ and the unbounded region containing $g_4$. Let us call the unbounded region as $R$. 
	Notice that, the bounded region is further subdivided into three regions , say $R_1, R_2$, and $R_3$ bounded by the cycles 
	$xg_1h_1g_2x$, $xg_2h_3g_3x$, and $g_1h_1g_2h_3g_3hg_1$, respectively.

	Notice that, any good vertex $g$ belonging to 
	$R_3$ must see $g_4$ via $h$. 
	However, then $x$ and $h$ will have three common good neighbors 
	$g_1, g_3, g_4$, creating a $\mathbb{F}_3$. 
	Therefore, $R_3$ cannot contain any good vertices.  
	
	Furthermore, any good vertex belonging to $R_1$ or $R_2$ must see $g_4$ via $x$, and any vertex belonging to $R$ must see $g_2$ via $x$. Thus, every good vertex, maybe except $h$ and $x$ itself, is adjacent to $x$. 
	As $x$ can have at most $2(2n+m)^2$ good neighbors due to Lemma~\ref{lem 2p alpha max}, there can be at most  $2(2n+m)^2+2$ good vertices in $H$, counting $x$ and $h$, a contradiction.  
\end{proof}

With the above lemma, we can further restrict the number of good neighbors of any vertex in $H$. 
\begin{lemma}\label{lem good neighbors}
	A vertex $x \in V(H)$ can have at most $2(2n+m)$ good neighbors in $H$.
\end{lemma}

\begin{proof}
Suppose $x$ has at least $2(2n+m)+1$ good neighbors. Note that $2(2n+m)+1 \geq 5$ as $2n+m \geq 2$.  
By the Pigeonhole principle, $x$ must have at least three  $\alpha$-neighbors for some adjacency type $\alpha$. 
According to Lemma~\ref{lem at most 3 alpha+2beta}, this is not possible, a contradiction. 
\end{proof}

% Finally, we are ready to prove Theorem~\ref{thm main}. 

%\bigskip

%\noindent \textit{Proof of Theorem~\ref{thm main}.}  

\begin{lemma}\label{lem help vertex}
The $(n,m)$-graph $H$ has no helper vertex of degree three. 
\end{lemma}

\begin{proof}
Suppose $u$ is a helper vertex in $H$ of degree three having neighbors $u_1, u_2$ and $u_3$. We delete the vertex $u$ and add three degree two helper vertices $h_{12}, h_{23}, h_{31}$ such that $h_{ij}$ is adjacent to $u_i$ and $u_j$. Moreover, if $u_i$ (resp., $u_j$) is an $\alpha$-neighbor of $u$ then it is also an $\alpha$-neighbor of $h_{ij}$. Let the so obtained graph be $H^*$.

Note that $H^*$ is a triangle-free planar graph satisfying $n_3(H^*) < n_3(H)$. 
Observe that, 
if two good vertices of $H$ see each other directly or via a helper other than $u$, then they see each other in the same way in $H^*$.
On the other hand, the only vertices that see each other via $u$ in $H$ are $u_i$ and $u_j$ for $i,j \in \{1,2,3\}$ while in $H^*$ they see each other via $h_{ij}$.  
Hence, $\omega_{r(n,m)}(H^*) \geq  \omega_{r(n,m)}(H)$ as the maximum $(n,m)$-relative clique $R$ in $H$ remains an $(n,m)$-relative clique in $H^*$ as well.  
  This contradicts the criticality of $H$. 
\end{proof} 

We recall a useful result for planar graphs.

\begin{theorem}[\cite{jendrol2013light}]\label{thm lightedge}
	Every planar graph $G$ with $\delta(G) \geq 4$, contains an edge $xy$ with $\deg(x) = 4$ and $4 \leq \deg(y) \leq 7$, or contains an edge $xy$, with $\deg(x) = 5$ and $ 5 \leq \deg(y) \leq 6$. The bounds are the best possible. 
	
\end{theorem}

\begin{lemma}\label{lem good vertex}
	There exists a good vertex in $H$ of degree at most seven. Moreover, if all good vertices have a degree of at least seven, then there exists a  good vertex of degree seven in $H$ with a neighbor which is a helper vertex of degree four.
\end{lemma}

\begin{proof}
	 We know that there exists no helper vertex of degree zero or one in $H$ due to its criticality and degree three in $H$ by Lemma \ref{lem help vertex}. Let us delete each helper vertex of degree two in $H$ and put an edge between its neighbors. Let the so obtained graph be $H'$,  Note that, two distinct helper vertices of degree two in $H$ do not have the same common neighborhood as $H$ is critical. Thus $H'$ is a simple planar graph with helper vertices of minimum degree four. 
	 
	 Therefore, either $H$ contains a good vertex of degree at most three or $\delta(H) \geq 4$. If the former happens, then we are done. If the latter happens, then by Theorem \ref{thm lightedge}, $H$ contains an edge $xy$ such that $\deg(x) + \deg(y) \leq 11$. Without loss of generality, let $\deg(x) \leq \deg(y)$. If $x$ is a good vertex, we are done. If $x$ is a helper vertex, then $y$ has to be a good vertex as helper vertices are independent by Lemma \ref{lem Sind}. 
	 
	 Note that, in case every good vertex in $H'$ has a degree at least seven, by the above argument, there must exist a good vertex $y$ of degree seven with a neighbor $x$ which is a helper vertex of degree four. 
\end{proof}

\subsection{Analysis of the structure of $H$}
Let $v$ be the good vertex in $H$ having a degree at most seven. Such a vertex exists due to Lemma \ref{lem good vertex}. 
 Let $v_1, v_2, v_3, \cdots ,v_k$ be its neighbors where $k \leq 7$. Let $A_i$ be the set of good neighbors of $v_i$ other than $v$ such that none of the good vertices in $A_i$ is adjacent to $v_j$ for any $j < i$. Let us assume that the indexing of the good neighbors of $v$ is such that the $k$-tuple $(|A_1|, |A_2|, |A_3|, \cdots ,|A_k|)$ attains a maximum in the lexicographic ordering. By Lemma~\ref{lem good neighbors} $|A_i| \leq 2(2n+m)-1$ for all $i \in \{1,2, \cdots ,k\}$. We would like to improve it by showing that,
 for $ 1 \leq p < q < r \leq 7$, we have $|A_p| \cup |A_q| \cup |A_r| \leq 2(2n+m)$.

\begin{lemma}\label{lem improvement}
	For $ 1 \leq p < q < r \leq 7$, we have $|A_p| \cup |A_q| \cup |A_r| \leq 2(2n+m)$. 
\end{lemma}

\begin{proof}

Suppose on the contrary, for some $1 \leq p < q < r \leq 7$, $|A_p| + |A_q| + |A_r| \geq 2(2n+m) + 1$. Without loss of generality, let us assume that the neighbors $v_p, v_q$ and $v_r$ be arranged in an anti-clockwise manner around $v$. And let the associated set of neighbors of each $v_p, v_q$ and $v_r$ be $A_p = \{ p_1, p_2, \cdots, p_a \}, A_q = \{ q_1, q_2, \cdots, q_b \}$ and  $A_r = \{ r_1, r_2, \cdots, r_c\}$, respectively. Furthermore, let the vertices in each of these sets $A_i$ to be arranged in an anti-clockwise manner around $v_i$ in increasing order of indices under the fixed embedding of $H$ for all $i \in \{p,q,r\}$. Following the above construction, by our lexicographic ordering and by Lemma~\ref{lem good neighbors}, we have, $2(2n+m)-1 \geq |A_p| \geq |A_q| \geq |A_r|$. Also $|A_p| + |A_q| + |A_r| \geq 2(2n+m) + 1 \geq 7$  as $(2n+m) \geq 3$. Hence, $|A_p| \geq 3$ and  $|A_q| + |A_r| \geq 2$. As $|A_q| + |A_r| \geq 2$, either $|A_q| \geq |A_r| \geq 1$ or $|A_q| \geq 2$ and $|A_r| \geq 0$. We analyse the cases based on how $p$ in $A_p$ sees $q_1$ in $A_q$, where 
$p = p_{1}$ if $|A_p| = a = 3$, and $p=p_{a-3}$, if $|A_p| = a \geq 4$. See Fig.~\ref{fig config} for a pictorial representations of the two scenarios. 

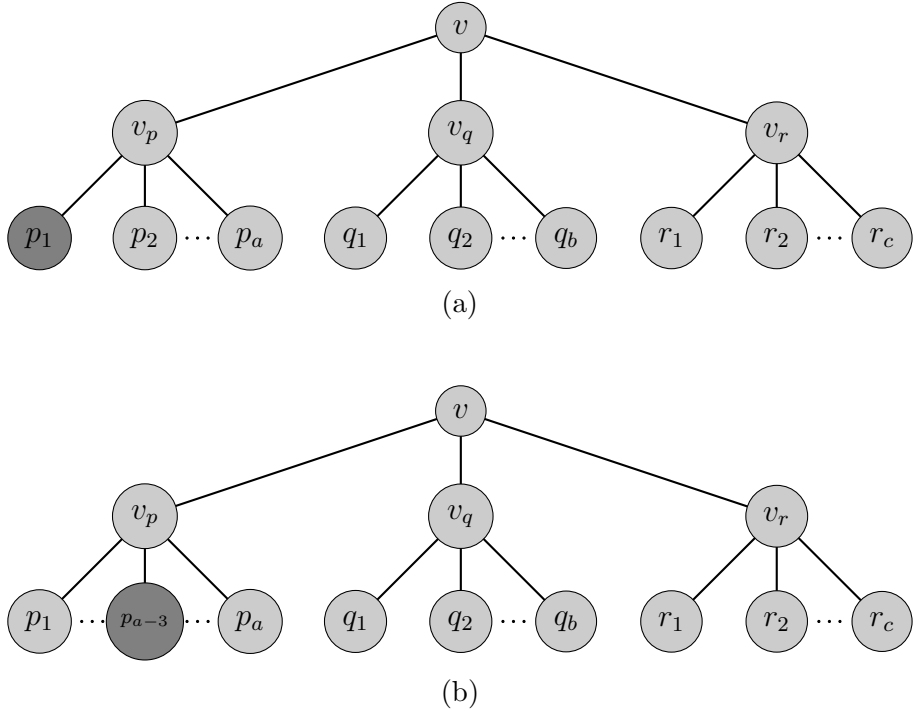
\begin{figure}
	\centering
             \begin{subfigure}[b]{1.05\textwidth}
             \centering
		\begin{tikzpicture}  
			[scale=.7,auto=center] % here, node/.style is the style pre-defined, that will be the default layout of all the nodes. You can also create different forms for different nodes.  

			\draw[thick] (6,4) -- (0,2);
			\draw[thick] (6,4) -- (6,2);
			\draw[thick] (6,4) -- (12,2);
			\draw[thick] (0,2) -- (-2,0);
			\draw[thick] (0,2) -- (0,0);
			\draw[thick] (0,2) -- (2,0);
			\draw[thick] (6,2) -- (4,0);
			\draw[thick] (6,2) -- (6,0);
			\draw[thick] (6,2) -- (8,0);
			\draw[thick] (12,2) -- (10,0);
			\draw[thick] (12,2) -- (12,0);
			\draw[thick] (12,2) -- (14,0);

			\node[circle,draw,fill=black!20] at (6,4) {$v$};
			\node[circle,draw,fill=black!20] at (0,2) {$v_p$};
				\node[circle,draw,fill=gray] at (-2,0) {$p_1$};
				\node[circle,draw,fill=black!20] at (0,0) {$p_2$};
				\draw[fill] (0.8,0) circle[radius=0.02];
				\draw[fill] (1,0) circle[radius=0.02];
				\draw[fill] (1.2,0) circle[radius=0.02];
				\node[circle,draw,fill=black!20] at (2,0) {$p_a$};
			\node[circle,draw,fill=black!20] at (6,2) {$v_q$};
				\node[circle,draw,fill=black!20] at (4,0) {$q_1$};
				\node[circle,draw,fill=black!20] at (6,0) {$q_2$};
				\draw[fill] (6.8,0) circle[radius=0.02];
				\draw[fill] (7,0) circle[radius=0.02];
				\draw[fill] (7.2,0) circle[radius=0.02];
				\node[circle,draw,fill=black!20] at (8,0) {$q_b$};
			\node[circle,draw,fill=black!20] at (12,2) {$v_r$};
				\node[circle,draw,fill=black!20] at (10,0) {$r_1$};
				\node[circle,draw,fill=black!20] at (12,0) {$r_2$};
				\draw[fill] (12.8,0) circle[radius=0.02];
				\draw[fill] (13,0) circle[radius=0.02];
				\draw[fill] (13.2,0) circle[radius=0.02];
				\node[circle,draw,fill=black!20] at (14,0) {$r_c$};
	
		\end{tikzpicture} 
  \caption{  }
	\label{fig config}
  \end{subfigure}

  \vspace{0.8cm}

\begin{subfigure}[b]{1.05\textwidth}
		\centering
			
		\begin{tikzpicture}  
			[scale=.7,auto=center] % here, node/.style is the style pre-defined, that will be the default layout of all the nodes. You can also create different forms for different nodes.  

			\draw[thick] (6,4) -- (0,2);
			\draw[thick] (6,4) -- (6,2);
			\draw[thick] (6,4) -- (12,2);
			\draw[thick] (0,2) -- (-2,0);
			\draw[thick] (0,2) -- (0,0);
			\draw[thick] (0,2) -- (2,0);
			\draw[thick] (6,2) -- (4,0);
			\draw[thick] (6,2) -- (6,0);
			\draw[thick] (6,2) -- (8,0);
			\draw[thick] (12,2) -- (10,0);
			\draw[thick] (12,2) -- (12,0);
			\draw[thick] (12,2) -- (14,0);

			\node[circle,draw,fill=black!20] at (6,4) {$v$};
			\node[circle,draw,fill=black!20] at (0,2) {$v_p$};
				\node[circle,draw,fill=black!20] at (-2,0) {$p_1$};
				\draw[fill] (-1.2,0) circle[radius=0.02];
				\draw[fill] (-1,0) circle[radius=0.02];
				\draw[fill] (-0.8,0) circle[radius=0.02];
                    \node[circle,draw,fill=gray] at (0,0) {\scriptsize $p_{a-3}$};
                    \draw[fill] (0.8,0) circle[radius=0.02];
				\draw[fill] (1,0) circle[radius=0.02];
				\draw[fill] (1.2,0) circle[radius=0.02];
				\node[circle,draw,fill=black!20] at (2,0) {$p_a$};
			\node[circle,draw,fill=black!20] at (6,2) {$v_q$};
				\node[circle,draw,fill=black!20] at (4,0) {$q_1$};
				\node[circle,draw,fill=black!20] at (6,0) {$q_2$};
				\draw[fill] (6.8,0) circle[radius=0.02];
				\draw[fill] (7,0) circle[radius=0.02];
				\draw[fill] (7.2,0) circle[radius=0.02];
				\node[circle,draw,fill=black!20] at (8,0) {$q_b$};
			\node[circle,draw,fill=black!20] at (12,2) {$v_r$};
				\node[circle,draw,fill=black!20] at (10,0) {$r_1$};
				\node[circle,draw,fill=black!20] at (12,0) {$r_2$};
				\draw[fill] (12.8,0) circle[radius=0.02];
				\draw[fill] (13,0) circle[radius=0.02];
				\draw[fill] (13.2,0) circle[radius=0.02];
				\node[circle,draw,fill=black!20] at (14,0) {$r_c$};
	
		\end{tikzpicture} 
  \caption{  }
	\label{fig config}
  \end{subfigure}

 \caption{A pictorial representations of the two scenarios dealt with in the proof of Lemma~\ref{lem improvement}.}
 \label{fig main fig}
\end{figure}

\medskip 

\noindent \textit{Case 1: If $p$ sees $q_1$ by being adjacent to it.} Let the region bounded by the cycle $pv_pvv_qq_1p$ containing $p_{a-1}$ be $D_1$. Now, $q_2$ (resp., $r_1$) can see $p_i$ only via the vertices of $D_1$, where $i=\{a-1,a\}$. 
    However, $q_2$ (resp., $r_1$)  cannot see $p_i$ via $p$ or $v$ as otherwise a triangle will be formed. Also, $q_2$ (resp., $r_1$) cannot see $p_i$ via $v_p$ as $v_p$ cannot be adjacent to $q_2$ (resp., $r_1$) due to the definition of $A_p$. 
    Moreover, $q_2$ (if it exists) cannot see $p_i$ via $q_1$ as otherwise a triangle will be formed. 
    Thus, the only possibility for $q_2$ to see $p_{a-1}$ and $p_{a}$ is via $v_q$. If $q_2$ see both $p_{a-1}, p_a$ via $v_q$, then $v_p,v_q$ along with their three common good neighbors $v, p_{a-1}, p_a$  will form an $\mathbb{F}_3$, which is not possible. Furthermore, $r_1$ (if it exists) cannot see $p_i$ via $v_q$ due to the definition of $A_q$. 
    Thus, the only possibility for $r_1$ to see $p_{a-1}$ and $p_{a}$ is via $q_1$. If $r_1$ sees both $p_{a-1}, p_{a}$  via $q_1$, then $v_p,q_1$ along with their three common good neighbors $p,p_{a-1},p_a$  will form an $\mathbb{F}_3$, which is not possible.  
    Therefore, it is not possible for $p$ to see $q_1$ by being adjacent to it.

  \medskip

\noindent \textit{Case 2: If $p$ sees $q_1$ via $v_q$.} 
Let the region bounded by the cycle $pv_pvv_qp$ containing $p_{a-1}$ be $D_2$. Note that $q_1$ can see $p_{a-1}$ only via the vertices of $D_2$. However, $q_1$ cannot see $p_{a-1}$ via $p$ or $v$  as otherwise a triangle will be formed. Also, $q_1$ cannot see $p_{a-1}$ via $v_p$ by the definition of $A_p$. Thus, the only possibility for $q_1$ to see $p_{a-1}$ is via $v_q$. If $q_1$ sees $p_{a-1}$ via $v_q$, then $v_p,v_q$ along with their three common good neighbors $v,p,p_{a-1}$  will form an $\mathbb{F}_3$, which is not possible.
    Therefore, it is not possible for $p$ to see $q_1$ via $v_q$.
    
\medskip  

\noindent \textit{Case 3: If $p$ sees $q_1$ via $v_r$.} 
Let the region bounded by the cycle $pv_pvv_qq_1v_rp$ containing $p_a$ be $D_3$. 
    Note that, $q_2$ (if it exists) can see $p_i$ only via the vertices of $D_3$ where $i \in \{a-1,a\}$. However, $q_2$ cannot see $p_i$ via $p, v$ or $q_1$ as otherwise a triangle  will be formed.  Also, $q_2$ cannot see $p_i$ via $v_p$ by the definition of $A_p$. Notice that $q_2$ cannot see $p_i$ via $v_r$, as otherwise $v_p,v_r$ along with their three common good neighbors $p,p_i,v$ will form an $\mathbb{F}_3$. 
    Thus, the only possibility for $q_2$ to see $p_i$ is via $v_q$. If both $p_{a-1}, p_{a}$, sees $q_2$ via $v_q$, then $v_p,v_q$ along with their three common good neighbors $p,p_{a-1},p_a$  will form an $\mathbb{F}_3$, which is not possible.
    Thus, $q_2$ does not exist. Hence, $|A_q| = 1$ and $r_1$ must exist. Observe that as $v_r$ is adjacent to $p$, $q_1$ and $r_1$, it contradicts the lexicographic maximality of     $(|A_1|, |A_2|, \cdots, |A_k|)$. 
    Therefore, it is not possible for $p$ to see $q_1$ via $v_r$.

\medskip  

\noindent \textit{Case 4: If $p$ sees $q_1$ via a vertex $h \notin \{v,v_p,v_q,v_r\} \cup A_p \cup A_q \cup A_r$.} Let the region bounded by the cycle $pv_pvv_qq_1hp$ containing $p_a $ be $D_4$. Now $q_2$ (resp., $r_1$) can see $p_i$ only via the vertices of $D_4$, where $i \in \{a-1, a\}$ if $a=3$ and $i \in \{a-2, a-1, a\}$ if $a \geq 4$. However, note that $q_2$ (resp., $r_1$) cannot see $p_i$ via $p$ or $v$ as otherwise a triangle will be formed. Also, $q_2$ (resp., $r_1$) cannot see $p_i$ via $v_p$ by the definition of $A_p$. Moreover, $q_2$ cannot see $p_i$ via $q_1$ to avoid creating a triangle and $r_1$ cannot see $p_i$ via $v_q$ by the definition of $A_q$. 
     Thus, $q_2$ (if it exists), must see exactly one of $p_{a-1}, p_a$ via $v_q$ and the other via $h$ in order to avoid creating an $\mathbb{F}_3$. Similarly, $r_1$ (if it exists), must see exactly one of $p_{a-1}, p_a$ via $q_1$ and the other via $h$ in order to avoid creating an $\mathbb{F}_3$. Hence, if both $q_2$ and $r_1$ existed, a triangle will be created among the vertices $v_q, q_1$ and one of the $p_i$s. So, at most one of $q_2$ or $r_1$ can exist, which implies that $a \geq 4$. Thus, $q_2$ (resp., $r_1$) will see $p_{a-2}$ via $v_q$ (resp., $q_1$) or $h$, creating an $\mathbb{F}_3$ which is not possible. Hence, it is not possible for $p$ to see $q_1$ via $h$.

\medskip 

     \noindent \textit{Case 5: If $p$ sees $q_1$ via $r_j$ for some $j \in \{1,2,3,\cdots,c\}$.} Let the region bounded by the cycle $pv_pvv_rr_jp$ containing $p_{a-1}$ be $D_5$. Notice that this bounded region is subdivided into two regions: $(i)$ $D_{51}$ bounded by the cycle $pv_pvv_qq_1r_jp$ containing $p_{a-1}$ and $(ii)$ $D_{52}$ bounded by the cycle $q_1v_qvv_rr_jq_1$ not containing $p_{a-1}$.

     Suppose there exists some $r \in A_r$ distinct from $r_j$. If $r$ does not belong to the region $D_{52}$, and hence does not belong to $D_5$, then $r$ can see $p_i$ only via the vertices of $D_5$, where $i \in \{a-1,a\}$. 
     However, $r$ cannot see $p_i$ via $p, v$ or $r_j$, as otherwise a triangle will be formed. Also, $r$ cannot see $p_i$ via $v_p$ by the definition of $A_p$, and $v_r$ due to the planarity of $H$. Thus $r$ must be contained in the region $D_{52}$. Now, $r$ can see $p_i$ only via the vertices of  $D_{52}$. 
     However, we have seen that $r$ cannot see $p_i$ via $v,r_j$ or $v_r$. Notice that $r$ cannot see both $p_{a-1}$ and $p_{a}$ via $v_q$ as otherwise, a $\mathbb{F}_3$ will be formed. Also, $r$ cannot see one of $p_{a-1}$, $p_a$ via $q_1$ and the other one via $v_q$, as otherwise  $v_q,q_1,r$ will form a triangle. So the only possibility for $r$ to see $p_i$s is via $q_1$. 
     However, if $r$ sees $p_i$ via $q_1$, then $q_2$ (it exists as $|A_q| \geq |A_r| \geq 2$ in this case) which is in $D_{52}$ can see $p_i$ only via $q_1$, which is not possible, as otherwise, a triangle will be formed. 
     Thus, $|A_r| = 1$ and $r_j=r_1$.

    If $|A_q| \geq 3$, then $q_2$ and $q_3$ belongs to the region $D_{52}$. Note that, $q_j$ can see $p_i$ only via $v_q$ or $r_1$ as otherwise, either a triangle will be formed or planarity of $H$ will be violated, where $j = \{2,3\}$ and $i = \{a-1,a\}$. 
     Also, if both $p_{a-1}, p_a$ are adjacent to $v_q$ (resp., $r_1$) or if both $q_2, q_3$ are adjacent to $r_1$, an $\mathbb{F}_3$ will be formed. Thus, there is no possible way for the all the $q_j$s to see the all the $p_i$s. Hence, $|A_q| \leq 2$. 
     If $|A_q| = 2$, then $|A_p| \geq 4$. Thus, $q_2$ will see at least two of $p_{a-2}, p_{a-1}$ and $p_a$ via $v_q$ or $r_1$, forming an $\mathbb{F}_3$. Hence, $|A_q| = 1$.

    If $|A_q| = 1$, then $|A_p| \geq 2(2n+m)-1 \geq 5$. Also, we will have $p = p_{a-3}$. In this case, $p_{s}$ must see $q_1$ via $r_1$, as otherwise either a triangle will be created or the definition of $A_p$ will be violated, where $s \in \{1, 2, \cdots, a-4\}$. Observe that, this will create an $\mathbb{F}_3$ unless $(2n+m)=3$ and $a=5$. 
    Therefore, for the remainder of the proof let us assume that $(2n+m)=3$ and $a=5$.
    
    Observe that $q_1$ cannot see $p_i$ for $i = \{a-2,a-1,a\}$ via $p$ or $v$, as otherwise a triangle will be formed. 
    Also, $q_1$ cannot see $p_i$ via $v_p$ by the definition of $A_p$ or via $r_1$ as it will create an $\mathbb{F}_3$.  
    If $p_{a-2}$ see $q_1$ via $v_q$ or by being adjacent, then $p_{a-1}$ and $p_a$ cannot see $r_1$ without creating an $\mathbb{F}_3$. 
    Therefore, $p_{a-2}$ must see $q_1$ via some 
    $h \not\in \{v, v_p, v_q, v_r\} \cup A_p \cup A_q \cup A_r$. 
    This will force $p_{a-1}$ and $p_{a}$ to see $r_1$  via $q_1$. 
    
    Now we know that, apart from the good vertices among 
    $\{v, v_p, v_q, v_r\} \cup A_p \cup A_q \cup A_r$, we still have at least ten more good vertices. If such a good vertex belongs to the region $D_5$, then it has to be adjacent to $q_1$ or $r_1$ for seeing $p, p_1$ and $r_1$. If such a good vertex does not belong to the region $D_5$, then also it has to be adjacent to $r_1$ for seeing $q_1$. Thus, by the pigeonhole principle, $q_1$ or $r_1$ is adjacent to at least eight good vertices (recall, they are already having three good neighbors in form of $p_{a-1}, p_a, r_1$ or $p_1, p, q_1$, respectively). However, we know due to Lemma~\ref{lem good neighbors} that a vertex can have at most $2(2n+m)$ good neighbors, which in our case is six as $(2n+m)=3$. 
    Hence it is not possible for  $p$ to see $q_1$ via $r_j$ for some $j \in \{1,2,3,\cdots,c\}$.

\medskip  

From all the cases above, we have that $p$ cannot see $q_1$, which is a contradiction to our basic assumption. 
\end{proof}

\subsection{Concluding the proof}
We are now ready to prove the Theorem \ref{thm main}.

\medskip

\noindent \textit{Proof of Theorem \ref{thm main}.} Recall that $v$ is a good vertex in $H$ satisfying $k =\deg(v) \leq 7$. 

We count the number of possible good neighbors of $v$. If $k \leq 6$, then by Lemma~\ref{lem improvement} the total number of good vertices in the second neighborhood of $v$ is at most $ 4 (2n+m)$. Moreover, by Lemma \ref{lem good vertex}, the total number of good vertices in the closed neighborhood of $v$ is at most seven. Hence, the total number of good vertices in $H$ is, $$ |R| \leq 4 (2n+m) + 7 < 2 (2n+m)^2 + 2$$ for all  $(2n+m) \geq 3$, a contradiction.

  If $k = 7$, then by Lemma~\ref{lem improvement}, $v$ has at most six good neighbors and at least one neighbor which is a helper of degree four. Thus the total number of good vertices in the second neighborhood of $v$ is at most $4(2n+m) + 3$. Hence, the total number of good vertices in $H$ is, $$ |R| \leq 4 (2n+m) + 3 + 7 < 2 (2n+m)^2 + 2$$ for all  $(2n+m) \geq 4$, a contradiction.

  \medskip

   Now, for the case when $(2n+m) = 3$ and when $k = 7$, we give a detailed counting argument as follows. Since $k = 7$, the set of good vertices in the second neighborhood of $v$ is $A_1, A_2, A_3, \cdots A_7$. As $(2n+m) = 3$, we have $|A_i| \leq 2(2n+m) - 1 = 5$ for all $i \in \{1,2, \cdots ,7\}$. Note that, $|A_1| \geq 3$ as by Lemma~\ref{lem improvement}, there exists a neighbor of $v$ which is a helper vertex of degree four.  Moreover as $|A_1| \geq |A_2| \geq \cdots \geq |A_7|$ and $|A_1| + |A_2| + |A_3| \leq 6$, we are left with the following possible cases to deal with. 

\begin{enumerate}[(i)]
	\item If $|A_1| = 5$, $|A_2| = 1$, $|A_i| = 0$ for all $i \geq 3$, then the total number of good vertices are at most $|R \cap N[v]| + \displaystyle\sum_{i = 1}^{7} |A_i| \leq  7 + 6 = 13$. 
	
	\item If $|A_1| = 4$, $|A_2| = 2$, $|A_i| = 0$ for all $i \geq 3$, then the total number of good vertices are at most $|R \cap N[v]| + \displaystyle\sum_{i = 1}^{7} |A_i| \leq  7 + 6 = 13$. 
		
	\item If $|A_1| = 4$ and $|A_i| \leq 1$, for all $i \geq 2$, then the total number of good vertices are at most $|R \cap N[v]| + \displaystyle\sum_{i = 1}^{7} |A_i| \leq  7 + 10 = 17$. 
			
	\item If $|A_1| = 3$, $|A_2| = 3$, $|A_i| = 0$ for all $i \geq 3$, then the total number of good vertices are at most $|R \cap N[v]| + \displaystyle\sum_{i = 1}^{7} |A_i| \leq  7 + 6 = 13$. 
				
	\item If $|A_1| = 3$, $|A_2| \leq 2$, $|A_i| \leq 1$ for all $i \geq 3$, then the total number of good vertices are at most $|R \cap N[v]| + \displaystyle\sum_{i = 1}^{7} |A_i| \leq  7 + 10 = 17$. 				
	\end{enumerate}
	
%	\item[(ii)] If $|A_1| < 5$, then for  $|A_1| = 4$ and $|A_i| \leq 1$ for all other $i's$,  the total number of good vertices are at most $ 4 + 6 + 7 + 1 = 18$, 
%	\item[(iii)] If $|A_1| = 3$, and let $|A_2| = 2$ and $|A_i| \leq 1$ for other $i's$, then the total number of good vertices are at most $ 3 + 2 + 5 + 7 + 1 = 18$,
%	\item[(iv)] If $|A_1| = 3$, then $|A_i| = 1$ for all $i's$, then we have the total number of good vertices are at most $ 3 + 6 + 7 + 1 = 17$, 
%	\item[(v)] If $|A_1| = 2$ and if $|A_i| \leq 1$ for all $i's$, then we have the total number of good vertices are at most $ 8 + 7 + 1 = 16$, 
%	\item[(vi)] If $|A_1| = 2$ and if $|A_i| \leq 2$ for all $ i's$,  then we have, $14 + 7 + 1 = 22$ but if all $A_i$'s have two good vertices in particular, the vertex $h$ has $3$ good vertices including $v$ which is not possible in $H$ as $H$ has no degree $3$ helper, due to Lemma \ref{lem help vertex}.

In all the cases above, we see that the number of good vertices is strictly less than 
$2(2n+m)^2 + 2 = 20$ since $(2n+m )= 3$, a contradiction. \qed

\section{Conclusions}\label{sec conclusions}
In the following, we present a table which describes what is known regarding the $(n,m)$-relative and $(n,m)$-absolute clique number of the families $\mathcal{P}_g$ of planar graphs having 
girth at least $g$. The previous works had successfully found the exact values of $\omega_a{(n,m)}(\mathcal{P}_g)$ for all $g \geq 4$ and 
$\omega_r{(n,m)}(\mathcal{P}_g)$ for all $g \geq 5$. Therefore, finding the exact values were left open only for
$\omega_a{(n,m)}(\mathcal{P}_3)$, 
$\omega_r{(n,m)}(\mathcal{P}_3)$ 
and $\omega_r{(n,m)}(\mathcal{P}_4)$.
In this article, we find the exact value  for 
$\omega_r{(n,m)}(\mathcal{P}_4)$.

\begin{table}[h!]
\centering
\begin{tabular}{|c|c|c|}
\hline
    $\mathcal{P}_g$, for $g$  &
    $\omega_{a(n,m)}(\mathcal{P}_g)$ & $\omega_r{(n,m)}(\mathcal{P}_g)$  \\
    \hline
    $=3$  & \textcolor{blue}{$\left[ 3p^2 + p + 1, 9p^2 + 2p + 2 \right]$}~\cite{bensmail2017analogues} & \textcolor{blue}{$\left[ 3p^2 + p + 1, 42p^2 - 11 \right]$}~\cite{chakraborty2021clique}  \\
    \hline
    $=4$  &   $p^2 + 2$~\cite{chakraborty2021clique} & \textbf{\textcolor{red}{$2p^2+2$}} (Theorem~\ref{thm main}) \\
    \hline
    $=5$   & $\max{(p+1, 5)}$~\cite{chakraborty2021clique} & $\max{(p+1, 6)}$~\cite{chakraborty2021clique} \\
    \hline
    $=6$  &  $p + 1$~\cite{chakraborty2021clique} & $\max{(p+1, 4)}$~\cite{chakraborty2021clique} \\
    \hline
    $ \geq 7 $ & $p+1$~\cite{chakraborty2021clique}& $p+1$~\cite{chakraborty2021clique} \\
    \hline
\end{tabular}

\caption{This is the list of all known lower and upper bounds for 
    $\omega_{a(n,m)}(\mathcal{P}_g)$, $\omega_r{(n,m)}(\mathcal{P}_g)$ where $\mathcal{P}_g$ denotes the family of planar graphs having girth at least $g$. Moreover, the list captures the bounds for all $(n,m) \neq (0,1)$ where $p=2n+m$.}
\label{table results}
\end{table}

We know that 
$$2(2n+m)^2 \leq \chi_{n,m}(\mathcal{P}_4) \leq 5(2n+m)^4$$
where the lower bound is by Theorem~\ref{thm chi} and the upper bound follows from the upper bound of the $(n,m)$-chromatic number of planar graphs due to Ne\v{s}et\v{r}il and Raspaud~\cite{nevsetvril2000colored}. As the $(n,m)$-relative clique number of the triangle-free planar graphs is 
quadratic in $(2n+m)$, the following may be an interesting question to ask:

\begin{question}
Does there exist a constant $c$ such that 
$\chi_{n,m}(\mathcal{P}_4) \leq c(2n+m)^2$ for all $2n+m \geq 3$?
\end{question}

\bigskip

\noindent \textbf{Acknowledgements:} This work is partially supported by IFCAM project ``Applications of graph homomorphisms''
(MA/IFCAM/18/39), and SERB-MATRICS ``Oriented chromatic and clique number of planar graphs'' (MTR/2021/000858).

\medskip 

\noindent \textit{Note:} The authors have no relevant financial or non-financial interests to disclose. The manuscript has no associated data. All the authors have contributed equally and their names are listed in the alphabetical order of their last names.

\bibliographystyle{plain}
\bibliography{bibfile.bib}

\end{document}